\documentclass{preprint}

\usepackage{amsmath,amsfonts}
\usepackage{mhequ}
\usepackage{mhenvs}
\usepackage{times}
\usepackage{verbatim}

\newcommand{\CC}{{\mathcal C}}

\newcommand{\CD}{{\mathcal D}}
\newcommand{\CH}{{\mathcal H}}
\newcommand{\CL}{{\mathcal L}}
\newcommand{\CN}{{\mathcal N}}

\newcommand{\CR}{{\mathcal R}}
\renewcommand{\d}{\partial}
\newcommand{\e}{{\textrm e}}
\newcommand{\E} {\mathbb{E}}
\newcommand{\eg}{\textit{e.g.}\ }
\let\eps=\epsilon
\newcommand{\eref}[1]{(\ref{#1})}
\newcommand{\grad}{\nabla}
\newcommand{\ie}{\textit{i.e.}\ }
\DeclareMathOperator{\Law}{law}
\newcommand{\lhs}{\hskip1cm&\hskip-1cm}
\newcommand{\N}{{\mathbb N}}
\renewcommand{\P} {\mathbb{P}}

\newcommand{\R}{{\mathbb R}}
\newcommand{\scal}[1]{\langle#1\rangle}

\def\TV{{\mathrm{TV}}}

\numberwithin{equation}{section}

\begin{document}

\title{Analysis of SPDEs Arising in Path Sampling\\
Part I: The Gaussian Case}

\author{M. Hairer\inst{1}, A. M. Stuart\inst{1}, J. Voss \inst{1}\fnmsep\inst{2}, and P. Wiberg \inst{1}\fnmsep\inst{3}}

\institute{
Mathematics Institute, The University of Warwick, Coventry CV4 7AL, UK \and
Previously at University of Kaiserslautern, 67653 Kaiserslautern, Germany; 
supported by Marie Curie Fellowship HPMT-CT-2000-00076 \and
Now at Goldman-Sachs, London
}
\maketitle


\begin{abstract}
  In many applications it is important to be able to sample paths of
  SDEs conditional on observations of various kinds.  This paper
  studies SPDEs which solve such sampling problems.  The SPDE may
  be viewed as an infinite dimensional analogue of the Langevin SDE
  used in finite dimensional sampling.  Here the theory is developed
  for conditioned Gaussian processes for which the resulting SPDE is
  linear.  Applications include the Kalman-Bucy filter/smoother.  A
  companion paper studies the nonlinear case, building on the linear
  analysis provided here.
\end{abstract}


\section{Introduction}

An important basic concept in sampling is Langevin dynamics: suppose
a target density~$p$ on~$\R^d$ has the form $p(x) = c
\exp\bigl(-V(x)\bigr)$.  Then the stochastic differential equation~(SDE)
\begin{equ}[e:Langevin]
  \frac{dx}{dt} = -\grad V(x) + \sqrt{2} \,\frac{dW}{dt}
\end{equ}
has $p$ as its invariant density.  Thus, assuming that
\eref{e:Langevin} is ergodic, $x(t)$ produces samples from the target
density $p$ as $t \to \infty$.  (For details see, for example,
\cite{Robert-Casella99}.)

In~\cite{Stuart-Voss-Wiberg04} we give an heuristic approach to
generalising the Langevin method to an infinite dimensional setting.  We
derive stochastic partial differential equations (SPDEs) which are the
infinite dimensional analogue of~\eref{e:Langevin}. These SPDEs sample
from paths of stochastic differential equations, conditional on
observations.  Observations which can be incorporated into this
framework include knowledge of the solution at two points (bridges)
and a set-up which includes the Kalman-Bucy filter/smoother.  For
bridge sampling the SPDEs are also derived
in~\cite{Reznikoff-VandenEijnden05}, their motivation being to understand
the invariant measures of SPDEs through bridge processes.

In the current paper we give a rigorous treatment of this SPDE based
sampling method when the processes to be sampled are linear and Gaussian.  The
resulting SPDEs for the sampling are also linear and Gaussian in this case. 
We find it useful to present the Gaussian theory of SPDE based sampling
for conditioned diffusions in a self-contained fashion for the
following reasons.

\begin{itemize}

\item For nonlinear problems the SPDE based samplers can be quite
  competitive.  A companion article~\cite{Hairer-Stuart-VossII} will
  build on the analysis in this paper to analyse SPDEs which sample
  paths from nonlinear SDEs, conditional on observations.  The
  mathematical techniques are quite different from the Gaussian
  methods used here and hence we present them in a separate paper.
  However the desired path-space measures there will be characterised
  by calculating the density with respect to the Gaussian measures
  calculated here.

\item We derive an explicit description of the Kalman/Bucy smoother
  via the solution of a linear two-point boundary value problem.  This
  is not something that we have found in the existing literature; it
  is strongly suggestive that for off-line smoothing of Gaussian
  processes there is the potential for application of a range of fast
  techniques available in the computational mathematics literature, and
  different from the usual forward/backward implementation of the
  filter/smoother. See section \ref{S:KB}.

\item For Gaussian processes, the SPDEs studied here will not usually
  constitute the optimal way to sample, because of the time
  correlation inherent in the SPDE; better methods can be developed to
  generate independent samples by factorising the covariance operator.
  However these better methods can be viewed as a particular
  discretisation of the SPDEs written down in this paper, and this
  connection is of both theoretical interest and practical use,
  including as the basis for algorithms in the nonlinear case. See
  section \ref{S:disc} and \cite{Roberts-Stuart-Voss}.

\end{itemize}

In section~\ref{S:gaussian} of this article we will develop a general
MCMC method to sample from a given Gaussian process.  It transpires
that the distribution of a centred Gaussian process coincides with the
invariant distribution of the $L^2$-valued SDE
\begin{equ}
\label{e:L2SDE}
\frac{dx}{dt}
  = \CL x -\CL m + \sqrt{2} \,\frac{dw}{dt}
  \qquad \forall t \in (0,\infty),
\end{equ}
where $\CL$ is the inverse of the covariance operator, $m$ is the mean
of the process and $w$ is a cylindrical Wiener process.  

The first sampling problems we consider are governed by paths
of the $\R^d$-valued linear SDE
\begin{equ}[e:SDE0]
  \frac{dX}{du}(u) = AX(u) + B\,\frac{dW}{du}(u)  \qquad \forall u \in [0,1]
\end{equ}
subject to observations of the initial point $X(0)$, as well as
possibly the end-point~$X(1)$.  Here we have $A,B\in\R^{d\times d}$
and $W$ is a standard $d$-dimensional Brownian motion.  Since the SDE
is linear, the solution~$X$ is a Gaussian process.  Section~3
identifies the operator~$\CL$ in the case where we sample solutions
of~\eref{e:SDE0}, subject to end-point conditions.  In fact, $\CL$ is
a second order differential operator with boundary conditions
reflecting the nature of the observations and thus we can
write~\eqref{e:L2SDE} as an~SPDE.

In section 4 we study the situation where two processes $X$ and $Y$
solve the linear system of SDEs
\begin{equs}[0]
  \frac{dX}{du}(u)
    = A_{11} X(u) + B_{11}\,\frac{dW_x}{du}(u) \\
  \frac{dY}{du}(u)
    = A_{21} X(u) + B_{22}\,\frac{dW_y}{du}(u)
\end{equs}
on $[0,1]$ and we want to sample paths from the distribution of $X$
(the signal) conditioned on~$Y$ (the observation).  Again, we identify
the operator $\CL$ in~\eref{e:L2SDE} as a second order differential
operator and derive an SPDE with this distribution as its invariant
distribution.  We also give a separate proof that the mean of the
invariant measure of the SPDE coincides with the standard algorithmic
implementation of the Kalman-Bucy filter/smoother through
forward/back\-ward sweeps.

Section \ref{S:disc} contains some brief remarks concerning the
process of discretising SPDEs to create samplers, and section
\ref{S:conc} contains our conclusions.

\medskip

To avoid confusion we use the following naming convention.  Solutions
to SDEs like~\eqref{e:SDE0} which give our target distributions are
denoted by upper case letters.  Solutions to infinite dimensional
Langevin equations like~\eqref{e:L2SDE} which we use to sample from
these target distributions are denoted by lower case letters.


\section{Gaussian Processes}
\label{S:gaussian}

In this section we will derive a Hilbert space valued~SDE to sample
from arbitrary Gaussian processes.

\medskip

Recall that a random variable $X$ taking values in a separable Hilbert
space $\CH$ is said to be \emph{Gaussian} if the law of $\scal{y,X}$
is Gaussian for every $y \in \CH$ (Dirac measures are considered as
Gaussian for this purpose).  It is called \emph{centred} if $\E
\scal{y,X} = 0$ for every $y \in \CH$.  Gaussian random variables are
determined by their mean $m = \E X \in \CH$ and their covariance
operator $\CC\colon\CH \to \CH$ defined by
\begin{equ}
 \scal{y,\CC x} = \E \bigl( \scal{y,X-m}\scal{X-m, x} \bigr).
\end{equ}
For details see \eg \cite[section~2.3.2]{DaPrato-Zabczyk92}.  The
following lemma (see \cite[proposition~2.15]{DaPrato-Zabczyk92})
characterises the covariance operators of Gaussian measures.

\begin{lemma} \label{L:covariance}
  Let $X$ be a Gaussian random variable on a separable Hilbert space.
  Then the covariance operator $\CC$ of~$X$ is self-adjoint, positive
  and trace class.
\end{lemma}

A Gaussian random variable is said to be \emph{non-degenerate} if
$\scal{y,\CC y} > 0$ for every $y \in \CH\setminus\{0\}$.  An
equivalent characterisation is that the law of $\scal{y,X}$ is a
proper Gaussian measure (\ie not a Dirac measure) for every $y \in
\CH\setminus\{0\}$.  Here we will always consider non-degenerate
Gaussian measures.  Then $\CC$ is strictly positive definite and we
can define~$\CL$ to be the inverse of~$-\CC$.  Since $\CC$ is trace
class, it is also bounded and thus the spectrum of $\CL$ is bounded
away from~$0$.

We now construct an infinite dimensional process which, in
equilibrium, samples from a prescribed Gaussian measure.  Denote by
$w$ the cylindrical Wiener process on~$\CH$.  Then one has formally
\begin{equ}[e:defW]
  w(t) = \sum_{n=1}^\infty \beta_n(t) \phi_n \quad \forall t \in (0,\infty),
\end{equ}
where for $n\in\N$ the $\beta_n$ are i.i.d.\ standard Brownian motions
and $\phi_n$ are the (orthonormal) eigenvectors of~$\CC$.  Note that
the sum~\eref{e:defW} does \textit{not} converge in $\CH$ but that one
can make sense of it by embedding $\CH$ into a larger Hilbert space in
such a way that the embedding is Hilbert-Schmidt. The choice of this
larger space does not affect any of the subsequent expressions (see
also \cite{DaPrato-Zabczyk92} for further details).

Given $\CC$ and $\CL$ as above, consider the $\CH$-valued SDE given by
\eqref{e:L2SDE}, interpreted in the following way:
\begin{equ}[e:defmild]
x(t) = m + \e^{\CL t} \bigl(x(0) - m\bigr)
                + \sqrt 2\int_0^t \e^{\CL(t-s)} \,dw(s).
\end{equ}
If $x\in C\bigl([0,T],\CH\bigr)$ satisfies~\eref{e:defmild} it is
called a {\em mild} solution of the SDE~\eref{e:L2SDE}.  We have the
following result.

\begin{lemma}\label{convergence}
  Let $\CC$ be the covariance operator and $m$ the mean of a
  non-degenerate Gaussian random variable $X$ on a separable Hilbert
  space~$\CH$.  Then the corresponding evolution
  equation~\eref{e:L2SDE} with $\CL=-\CC^{-1}$ has continuous
  $\CH$-valued mild solutions.  Furthermore, it has a unique invariant
  measure $\mu$ on $\CH$ which is Gaussian with mean~$m$ and
  covariance~$\CC$ and there exists a constant $K$ such that for every
  initial condition $x_0 \in \CH$ one has
  \begin{equ}
    \bigl\|\Law\bigl(x(t)\bigr) - \mu\bigr\|_\TV
      \le K \, \bigl(1+\|x_0 - m\|_\CH\bigr)
                 \exp\bigl({-\|\CC\|^{-1}_{\CH\to\CH} t}\bigr) ,
  \end{equ}
  where $\|\cdot\|_\TV$ denotes the total variation distance between
  measures.
\end{lemma}

\begin{proof}
The existence of a continuous $\CH$-valued solution of the SDE~\eref{e:L2SDE}
is established in \cite{Iscoe-Marcus-McDonald-Talagrand-Zinn90}. The
uniqueness of the invariant measure and the convergence rate in the
total variation distance follow by combining Theorems~6.3.3 and~7.1.1
from \cite{DaPrato-Zabczyk96}.  The characterisation of the invariant
measure is established in \cite[Thm~6.2.1]{DaPrato-Zabczyk96}.
\end{proof}

We can both characterise the invariant measure, and explain the
exponential rate of convergence to it, by using the Karhunen-Lo\`eve 
expansion. In particular we give an heuristic argument
which illustrates why Lemma~\ref{convergence} holds in the case $m=0$:
denote by $(\phi_n)_{n\in\N}$ an orthonormal basis of eigenvectors
of~$\CC$ and by $(\lambda_n)_{n\in\N}$ the corresponding eigenvalues.
If $X$ is centred it is possible to expand $X$ as
\begin{equ}[e:KL]
  X = \sum_{n=1}^{\infty} \alpha_n \sqrt{\lambda_n} \phi_n,
\end{equ}
for some real-valued random variables $\alpha_n$.  (In contrast to the
situation in~\eqref{e:defW} the convergence in \eref{e:KL} actually
holds in $L^2(\Omega, \P, \CH)$, where $(\Omega, \P)$ is the
underlying probability space.)  A simple calculation shows that the
coefficients $\alpha_n$ are i.i.d.\ $\CN(0,1)$ distributed random
variables.  The expansion~\eref{e:KL} is called the Karhunen-Lo\`eve
expansion.  Details about this construction can be found
in~\cite{Adler}.

Now express the solution~$x$ of~\eqref{e:L2SDE} in the basis
$(\phi_n)$ as
\begin{equ}
  x(t) = \sum_{n=1}^\infty \gamma_n(t) \phi_n.
\end{equ}
Then a formal calculation using \eqref{e:defW} and~\eqref{e:L2SDE}
leads to the SDE
\begin{equ}
  \frac{d\gamma_n}{dt} = - \frac{1}{\lambda_n} + \sqrt{2} \frac{d\beta_n}{dt}
\end{equ}
for the time evolution of the coefficients~$\gamma_n$ and hence
$\gamma_n$ is ergodic with stationary distribution~$\CN(0,\lambda_n)$
for every~$n\in\N$.  Thus the stationary distribution
of~\eqref{e:L2SDE} has the same Karhunen-Lo\`eve expansion as the
distribution of~$X$ and the two distributions are the same.

\bigskip

In this article, the Hilbert space $\CH$ will always be the space
$L^2\bigl([0,1], \R^d\bigr)$ of square integrable $\R^d$-valued
functions and the Gaussian measures we consider will be distributions
of Gaussian processes.  In this case the operator $\CC$ has a kernel
$C\colon [0,1]^2 \to \R^{d \times d}$ such that
\begin{equ}[e:defC]
(\CC x)(u) = \int_0^1 C(u,v) \, x(v)\,dv.
\end{equ}
If the covariance function $C$ is H\"older continuous, then the
Kolmogorov continuity criterion (see \eg
\cite[Thm~3.3]{DaPrato-Zabczyk92}) ensures that $X$ is almost surely a
continuous function from $[0,1]$ to $\R^d$.  In this case $C$ is given
by the formula
\begin{equ}
C(u,v) = \E \Bigl( \bigl(X(u)-m(u)\bigr)\bigl(X(v)-m(v)\bigr)^* \Bigr)
\end{equ}
and the convergence of the expansion~\eref{e:KL} is uniform with
probability one.

\begin{remark}\label{rem:conv}
  The solution of~\eref{e:L2SDE} may be viewed as the basis for 
an MCMC me\-thod for
  sampling from a given Gaussian process.  The key to exploiting this
  fact is the identification of the operator~$\CL$ for a given
  Gaussian process.  In the next section we show that, for a variety
  of linear SDEs, $\CL$ is a second order differential operator and
  hence~\eref{e:L2SDE} is a stochastic partial differential equation.
  If $\CC$ has a H\"older continuous kernel~$C$, it follows
  from~\eref{e:defC} and the relation $\CC = (-\CL)^{-1}$ that it
  suffices to find a differential operator~$\CL$ such that $C(u,v)$ is
  the Green's function of~$-\CL$.
\end{remark}


\section{Conditioned Linear SDEs}
\label{S:conditioned SDEs}

In this section we apply our sampling technique from
section~\ref{S:gaussian} to Gaussian measures which are given as the
distributions of a number of conditioned linear SDEs.  We condition
on, in turn, a single known point (subsection \ref{ssec:3.1}), a
single point with Gaussian distribution (subsection \ref{ssec:3.2})
and finally a bridge between two points (subsection \ref{ssec:3.3}).

\medskip

Throughout we consider the $\R^d$-valued SDE
\begin{equ}[e:SDE]
  \frac{dX}{du}(u) = AX(u) + B\,\frac{dW}{du}(u),   \qquad \forall u\in[0,1],
\end{equ}
where $A,B\in\R^{d\times d}$ and $W$ is the standard $d$-dimensional
Brownian motion.  We assume that the matrix $BB^*$ is invertible.  We
associate to \eref{e:SDE} the second order differential operator $L$
formally given by
\begin{equ}[e:opL]
L = ({\d_u} + A^*)
          (BB^*)^{-1}({\d_u} - A).
\end{equ}
When equipped with homogeneous boundary conditions through its domain
of definition, we will denote the operator~\eref{e:opL} by~$\CL$.  We
will always consider boundary conditions of the general form $D_0x(0)
= 0$ and $D_1x(1) = 0$, where $D_i = A_i\d_u + b_i$ are
first-order differential operators.

\begin{remark} \label{rem:BC}
  We will repeatedly write $\R^d$-valued SPDEs with inhomogeneous
  boundary conditions of the type
  \begin{equs}[0][e:formalBC]
    \d_t x(t,u) = L x(t,u) + g(u) + \sqrt{2} \,\d_t w(t,u)
                      \quad \forall (t,u) \in (0,\infty) \times [0,1], \\
    D_0 x(t,0) = a, \quad D_1 x(t,1) = b \qquad \forall t\in(0,\infty), \\
    x(0,u) = x_0(u) \qquad \forall u\in[0,1]
  \end{equs}
  where $g\colon [0,1]\to \R^d$ is a function, $\d_tw$ is space-time
  white noise, and $a,b\in\R^d$.  We call a process~$x$ a solution of
  this SPDE if it solves~\eref{e:defmild} with $x(0)=x_0$ where $\CL$ is
  $L$ equipped with the boundary conditions $D_0f(0) = 0$ and $D_1f(1)
  = 0$, and $m\colon[0,1]\to\R^d$ is the solution of the boundary
  value problem $-Lm=g$ with boundary conditions $D_0m(0)=a$ and
  $D_1m(1)=b$.

  To understand the connection between \eref{e:formalBC}
  and~\eref{e:defmild} note that, if $w$ is a smooth function, then
  the solutions of both equations coincide.
\end{remark}

\subsection{Fixed Left End-Point}
\label{ssec:3.1}

Consider the problem of sampling paths of~\eref{e:SDE} subject only
to the initial condition
\begin{equ}[e:c2lcond]
  X(0) = x^- \in \R^d.
\end{equ}
The solution of this SDE is a Gaussian process with mean
\begin{equ}[e:c2exp]
  m(u) = E\bigl(X(u)\bigr) = \e^{uA} x^-
\end{equ}
and covariance function
\begin{equ}[e:c1cov]
C_0(u,v)
  = \e^{uA} \Bigl( \int_0^{u\wedge v} \e^{-rA} B
                B^* \e^{-rA^*} \,dr \Bigr) \e^{vA^*}
\end{equ}
(see \eg \cite[section~5.6]{Karatzas-Shreve91} for reference).
Let~$\CL$ denote the differential operator~$L$ from~\eref{e:opL} with
the domain of definition
\begin{equ}[e:domL]
\CD(\CL) = \bigl\{ f\in H^2([0,1],\R^d) \bigm|
                 f(0) = 0, \frac{d}{du}f(1) = A f(1) \bigr\}.
\end{equ}

\begin{lemma} \label{Green}
With $\CL$ given by \eref{e:opL} and~\eref{e:domL} the function~$C_0$
is the Green's function for~$-\CL$.  That is
\begin{equ}
  L C_0(u,v) = - \delta(u-v)I
\end{equ}
and
\begin{equ}
  C_0(0,v) = 0, \quad \d_uC_0(1,v) = A\,C_0(1,v)
      \qquad \forall v\in(0,1).
\end{equ}
\end{lemma}

\begin{proof}
From~\eref{e:c1cov} it is clear that the left-hand boundary condition
$C_0(0,v)=0$ is satisfied for all $v \in [0,1]$. It also follows
that, for $u\neq v$, the kernel is
differentiable with derivative
\begin{equ}[e:oneThirdOfL]
\d_u C_0(u,v)
  = \begin{cases}
    A C_0(u,v) + B B^* \e^{-uA^*} \e^{vA^*},& \text{for $u<v$, and} \\
    A C_0(u,v) & \text{for $u>v$.}
  \end{cases}
\end{equ}
Thus the kernel~$C_0$ satisfies the boundary condition
$\d_u C_0(1,v)=AC_0(1,v)$ for all $v\in[0,1)$.

Equation~\eref{e:oneThirdOfL} shows
\begin{equ}[e:halfOfL]
(B B^*)^{-1}\bigl(\d_u - A\bigr) C_0(u,v)
  = \begin{cases}
    \e^{-uA^*} \e^{vA^*},& \text{for $u<v$, and} \\
    0 & \text{for $u>v$}
  \end{cases}
\end{equ}
and thus we get
\begin{equ}
LC_0(u,v)
  = \bigl(\d_u + A^*\bigr)
        (B B^*)^{-1}\bigl(\d_u - A\bigr) C_0(u,v)
  = 0  \qquad \forall u\neq v.
\end{equ}

Now let~$v\in(0,1)$.  Then we get
\begin{equ}
\lim_{u\uparrow v} \;
        (B B^*)^{-1}\bigl(\d_u - A\bigr) C_0(u,v)
  = I
\end{equ}
and
\begin{equ}
\lim_{u\downarrow v} \;
        (B B^*)^{-1}\bigl(\d_u - A\bigr) C_0(u,v)
  = 0
\end{equ}
This shows $LC_0(u,v) = -\delta(u-v)I$ for all $v\in(0,1)$.
\end{proof}

Now that we have identified the operator~$\CL=(-\CC)^{-1}$ we are in
the situation of~\lem{convergence} and can derive an SPDE to sample
paths of~\eref{e:SDE}, subject to the initial
condition~\eref{e:c2lcond}.  We formulate this result precisely in the
following theorem.

\begin{theorem} \label{samplingL}
For every $x_0\in\CH$ the $\R^d$-valued SPDE
\begin{subequations} \label{e:c2SPDE}
\begin{gather}
\d_t x(t,u) = L x(t,u) + \sqrt{2} \,\d_t w(t,u)
  \qquad \forall (t,u) \in (0,\infty)\times(0,1) \label{e:c2SPDEa} \\
x(t,0) = x^-, \quad \d_u x(t,1) = A x(t,1)
  \qquad \forall t\in(0,\infty) \label{e:c2SPDEb} \\
x(0,u) = x_0(u)
  \qquad \forall u\in[0,1]  \label{e:c2SPDEc}
\end{gather}
\end{subequations}
where $\d_t w$ is space-time white noise has a unique mild solution.
The SPDE is ergodic and in equilibrium samples paths of the
SDE~\eref{e:SDE} with initial condition~$X(0) = x^-$.
\end{theorem}

\begin{proof}
  The solution of SDE~\eref{e:SDE} with initial
  condition~\eref{e:c2lcond} is a Gaussian process where the mean~$m$
  is given by~\eref{e:c2exp}.  The mean~$m$ solves the boundary value
  problem~$Lm(u) = 0$ for all $u\in(0,1)$, $m(0) = x^-$ and $m'(1) =
  Am(1)$.  From \rem{rem:BC} we find that $x$ is a solution of the
  Hilbert space valued SDE~\eref{e:L2SDE} for this function~$m$.

  \lem{Green} shows that $\CL$, given by~\eref{e:opL} with the
  boundary conditions from~\eref{e:c2SPDEb}, is the inverse of~$-\CC$
  where $\CC$ is the covariance operator of the distribution we want
  to sample from (and with covariance function given
  by~\eref{e:c1cov}).  \lem{convergence} then shows that the
  SPDE~\eref{e:c2SPDE} is ergodic and that its stationary distribution
  coincides with the distribution of solutions of the SDE~\eref{e:SDE}
  with initial condition~$X(0) = x^-$.
\end{proof}

\subsection{Gaussian Left End-Point}
\label{ssec:3.2}

An argument similar to the one in section~\ref{ssec:3.1} deals with
sampling paths of~\eref{e:SDE} where $X(0)$ is a Gaussian random
variable distributed as
\begin{equ}[e:c3lcond]
  X(0) \sim \CN(x^-, \Sigma)
\end{equ}
with an invertible covariance matrix~$\Sigma\in\R^{d\times d}$ and
independent of the Brownian motion~$W$.

\begin{theorem} \label{samplingLn}
For every $x_0\in\CH$ the $\R^d$-valued SPDE
\begin{subequations} \label{e:c3SPDE}
\begin{gather}
\d_t x(t,u) = L x(t,u) + \sqrt{2} \,\d_t w(t,u)
  \qquad \forall (t,u) \in (0,\infty)\times(0,1) \label{e:c3SPDEa} \\
\d_u x(t,0) = A x(t,0) + BB^*\Sigma^{-1}(x-x^-),
  \quad \d_u x(t,1) = A x(t,1)
  \qquad \forall t\in(0,\infty) \label{e:c3SPDEb} \\
x(0,u) = x_0(u)
  \qquad \forall u\in[0,1] \label{e:c3SPDEc}
\end{gather}
\end{subequations}
where $\d_t w$ is space-time white noise has a unique mild solution.
The SPDE is ergodic and in equilibrium samples paths of the
SDE~\eref{e:SDE} with Gaussian initial condition~\eref{e:c3lcond}.
\end{theorem}

\begin{proof}
The solution~$X$ of SDE~\eref{e:SDE} with initial
condition~\eref{e:c3lcond} is a Gaussian process with
mean~\eref{e:c2exp} and covariance function
\begin{equ}[e:c3cov]
  C(u,v) = \e^{uA}\Sigma\e^{vA^*} + C_0(u,v),
\end{equ}
where $C_0$ is the covariance function from~\eref{e:c1cov} for the
case $X(0)=0$ (see Problem~6.1 in Section~5.6 of
\cite{Karatzas-Shreve91} for a reference).  The mean~$m$
from~\eref{e:c2exp} solves the boundary value problem~$Lm(u) = 0$ for
all $u\in(0,1)$ with boundary conditions $m'(0) = Am(0) +
BB^*\Sigma^{-1}(m(0)-x^-)$ and $m'(1) = Am(1)$.

In order to identify the inverse of the covariance operator~$\CC$ we
can use~\eref{e:oneThirdOfL} to find
\begin{equ}
\d_u C(u,v)
   = \begin{cases}
      A C(u,v) + B B^* \e^{-uA^*} \e^{vA^*},& \text{for $u<v$, and} \\
      A C(u,v) & \text{for $u>v$}
    \end{cases}
\end{equ}
and, since $C(0,v) = \Sigma\, \e^{v A^*}$, we get the boundary
conditions
\begin{equ}
  \d_u C(0,v) = A C(0,v) + B B^* \Sigma^{-1} C(0,v)
\end{equ}
and
\begin{equ}
  \d_u C(1,v) = A C(1,v).
\end{equ}
From $(\d_u-A) \e^{uA}\Sigma\e^{vA^*} = 0$ we also get
\begin{equ}
  LC(u,v) = L \e^{uA}\Sigma\e^{vA^*} + L C_0(u,v) = 0
\end{equ}
for all $u\neq v$ and $LC(u,v) = LC_0(u,v) = -\delta(u,v)I$ for all
$u,v\in(0,1)$.  Thus $C$ is again the Green's function for $-\CL$ and
the claim follows from \rem{rem:conv} and \lem{convergence}.
\end{proof}

\begin{remark}
  If $A$ is negative-definite symmetric, then the solution $X$ of
  SDE~\eref{e:SDE} has a stationary distribution which is a centred
  Gaussian measure with covariance $\Sigma = -\frac12 A^{-1}\,BB^*$.
  Choosing this distribution in~\eref{e:c3lcond}, the boundary
  condition~\eref{e:c3SPDEb} becomes
  \begin{equ}
    \d_u x(t,0) = - A x(t,0),
    \quad \d_u x(t,1) = A x(t,1)
    \qquad \forall t\in(0,\infty).
  \end{equ}
\end{remark}

\subsection{Bridge Sampling}
\label{ssec:3.3}

In this section we apply our sampling method to sample from solutions
of the linear SDE~\eref{e:SDE} with fixed end-points, \ie we sample
from the distribution of $X$ conditioned on
\begin{equ}[e:c4cond]
  X(0)=x^-, \quad X(1)=x^+.
\end{equ}
The conditional distribution transpires to be absolutely continuous
with respect to the Brownian bridge measure satisfying~\eref{e:c4cond}.

\medskip

Let $m$ and $C_0$ be the mean and covariance of the unconditioned
solution~$X$ of the SDE~\eref{e:SDE} with initial condition
$X(0)=x^-$.  As we will show in Lemma~\ref{gausscond} below, the
solution conditioned on $X(1)=x^+$ is again a Gaussian process.  The
mean and covariance of the conditioned process can be found by
conditioning the random variable $\bigl(X(u),X(v),X(1)\bigr)$ for
$u\leq v\leq 1$ on the value of $X(1)$.  Since this is a finite
dimensional Gaussian random variable, mean and covariance of the
conditional distribution can be explicitly calculated.  The result for
the mean is
\begin{equ}[e:c4exp]
  \tilde m(u) = m(u) + C_0(u,1)C_0(1,1)^{-1} \bigl(x^+-m(1)\bigr)
\end{equ}
and for the covariance function we get
\begin{equ}[e:c4cov]
  \tilde C(u,v) = C_0(u,v) - C_0(u,1) C_0(1,1)^{-1} C_0(1,v).
\end{equ}

\begin{theorem}
For every $x_0\in\CH$ the $\R^d$-valued SPDE
\begin{subequations} \label{e:c4SPDE}
\begin{gather}
\d_t x = L x + \sqrt{2} \,\d_t w
  \qquad \forall (t,u) \in (0,\infty)\times(0,1) \label{e:c4SPDEa} \\
x(t,0) = x^-, \quad x(t,1) = x^+
  \qquad \forall t\in(0,\infty) \label{e:c4SPDEb} \\
x(0,u) = x_0(u)
  \qquad \forall u\in[0,1]  \label{e:c4SPDEc}
\end{gather}
\end{subequations}
where $\d_t w$ is white noise has a unique mild solution.  The SPDE is
ergodic and in equilibrium samples paths of the SDE~\eref{e:SDE}
subject to the bridge conditions~\eref{e:c4cond}.
\end{theorem}

\begin{proof}
  The solution of the SDE~\eref{e:SDE} with boundary
  conditions~\eref{e:c4cond} is a Gaussian process where the mean
  $\tilde m$ is given by~\eref{e:c4exp} and the covariance function
  $\tilde C$ is given by~\eref{e:c4cov}.  From
  formula~\eref{e:halfOfL} we know $L C_0(u,1) = 0$ and thus $\tilde
  m$ satisfies $L\tilde m = L m = 0$.  Since $\tilde m(0) = x^-$ and
  $\tilde m(t) = m(1) + C_0(1,1) C_0(1,1)^{-1} \bigl(x^+ - m(1)\bigr)
  = x^+ $, the mean $\tilde m$ solves the boundary value problem
  $L\tilde m(u) = 0$ for all $u\in(0,1)$ with boundary conditions
  $\tilde m(0) = x^-$ and $\tilde m(1) = x^+$.

  It remains to show that $\tilde C$ is the Green's function for the
  operator~$L$ with homogeneous Dirichlet boundary conditions: we have
  $\tilde C(0,v) = 0$,
  \begin{equ}
    \tilde C(1,v) = C_0(1,v) - C_0(1,1)C_0(1,1)^{-1}C_0(1,v) = 0
  \end{equ}
  and using $L C_0(u,1) = 0$ we find
  \begin{equ}
    L \tilde C(u,v) = L C_0(u,v) = - \delta(u-v)I.
  \end{equ}
  This completes the proof.
\end{proof}


\goodbreak

\section{The Kalman-Bucy Filter/Smoother}
\label{S:KB}

Consider~\eref{e:SDE} with $X$ replaced by the $\R^m\times\R^n$-valued
process $(X,Y)$ and $A,B\in \R^{(m+n)\times (m+n)}$ chosen so as to obtain
the linear SDE
\minilab{e:kSDE}
\begin{equ}[e:kSDEa]
\frac{d}{du}\begin{pmatrix} X(u) \\ Y(u) \end{pmatrix}
  = \begin{pmatrix} A_{11} & 0 \\ A_{21} & 0 \end{pmatrix}
    \begin{pmatrix} X(u) \\ Y(u) \end{pmatrix}
  + \begin{pmatrix} B_{11} & 0 \\ 0 & B_{22} \end{pmatrix}
    \,\frac{d}{du}\!\begin{pmatrix} W_x(u) \\ W_y(u) \end{pmatrix}.
\end{equ}
We impose the conditions
\minilab{e:kSDE}
\begin{equ}[e:kSDEb]
  X_0 \sim \CN(x^-,\Lambda), \quad Y_0 = 0
\end{equ}
and try to sample from paths of $X$ given paths of~$Y$.  We derive an
SPDE whose invariant measure is the conditional distribution of $X$
given~$Y$.  Formally this SPDE is found by writing the SPDE for
sampling from the solution $(X,Y)$ of~\eref{e:kSDE} and considering
the equation for the evolution of~$x$, viewing $y\equiv Y$ as known.  This
leads to the following result.

\begin{theorem} \label{kalman-SPDE}
Given a path $Y$ sampled from~\eref{e:kSDE} consider the SPDE
\minilab{e:kSPDE}
\begin{equs}
\d_t x
  &= \Bigl((\d_u + A_{11}^*)
             (B_{11}B_{11}^*)^{-1}(\d_u - A_{11})\Bigr) x \\
  &\hskip20mm
   + A_{21}^*(B_{22}B_{22}^*)^{-1} \Bigl(\frac{dY}{du}-A_{21}x\Bigr)
   + \sqrt{2} \,\d_t w,
\label{e:kSPDEa}
\end{equs}
equipped with the inhomogeneous boundary conditions
\begin{equs} \minilab{e:kSPDE}
\d_u x (t,0)
      &= A_{11}x(t,0) + B_{11}B_{11}^* \Lambda^{-1} \bigl(x(t,0)-x^-\bigr), \\
  \d_u x(t,1) &= A_{11}x(t,1) \label{e:kSPDEb}
\end{equs}
and initial condition
\begin{equ} \minilab{e:kSPDE}
x(0,u) = x_0(u)
  \qquad \forall u\in[0,1].  \label{e:kSPDEc}
\end{equ}
Then for every $x_0\in\CH$ the SPDE has a unique mild solution and is
ergodic.  Its stationary distribution coincides with the conditional
distribution of $X$ given~$Y$ for $X,Y$ solving~\eref{e:kSDE}.
\end{theorem}

The proof of this theorem is based on the following three lemmas
concerning conditioned Gaussian processes.  After deriving these three
lemmas we give the proof of Theorem~\ref{kalman-SPDE}.  The section
finishes with a direct proof that the mean of the invariant measure
coincides with the standard algorithmic implementation of the
Kalman-Bucy filter/smoother through forward/backward sweeps
(this fact is implicit in Theorem~\ref{kalman-SPDE}).

\bigskip

\begin{lemma} \label{hilbert-schmidt} Let $\CH = \CH_1 \oplus \CH_2$
  be a separable Hilbert space with projectors $\Pi_i\colon
  \CH\to\CH_i$.  Let $\CC\colon \CH \to \CH$ be a positive
  definite, bounded, linear, self-adjoint operator and denote
  $\CC_{ij} = \Pi_i \CC \Pi_j^*$.  Then $\CC_{11} -
  \CC_{12}\CC_{22}^{-1}\CC_{21}$ is positive definite and if
  $\CC_{11}$ is trace class then the operator
  $\CC_{12}\CC_{22}^{-\frac12}$ is Hilbert-Schmidt.
\end{lemma}

\begin{proof}
Since $\CC$ is positive definite, one has
\begin{equ}
2 |\scal{\CC_{21} x, y}| \le \scal{x, \CC_{11}x} + \scal{y, \CC_{22}y},
\end{equ}
for every $(x,y)\in \CH$.  It follows that
\begin{equ}\label{e:basic}
|\scal{\CC_{21} x, y}|^2 \le \scal{x, \CC_{11}x} \scal{y, \CC_{22}y},
\end{equ}
and so
\begin{equ}\label{e:bb}
|\scal{\CC_{21} x, \CC_{22}^{-1/2} y}|^2 \le \scal{x, \CC_{11}x} \|y\|^2
\end{equ}
for every $y\neq 0$ in the range of $\CC_{22}^{1/2}$.
Equation~\eref{e:basic} implies that $\CC_{21} x$ is orthogonal to
$\ker \CC_{22}$ for every $x \in \CH_1$.  Therefore the operator
$\CC_{22}^{-1/2}\CC_{21}$ can be defined on all of $\CH_1$ and thus is
bounded.  Taking $y=\CC_{22}^{-1/2}\CC_{21}x$ in \eref{e:bb} gives
$\|\CC_{22}^{-1/2}\CC_{21}x\|^2 \leq \scal{x, \CC_{11}x}$ and thus
$\scal{x, (\CC_{11} - \CC_{12}\CC_{22}^{-1}\CC_{21}) x} \geq 0$ for every $x
\in \CH_1$.  This implies that $\CC_{22}^{-\frac12}\CC_{21}$
and $\CC_{12}\CC^{-\frac12}_{22}$ are both Hilbert-Schmidt, and
completes the proof.
\end{proof}

\begin{remark}
  Note that $\CC$ being strictly positive definite is not sufficient
  to imply that $\CC_{11} - \CC_{12}\CC_{22}^{-1}\CC_{21}$ is also
  strictly positive definite.  A counter-example can be constructed by
  considering the Wiener measure on $\CH = L^2([0,1])$ with $\CH_1$
  being the linear space spanned by the constant function $1$.
\end{remark}

\begin{lemma} \label{gausscond}
Let~$\CH=\CH_1\oplus\CH_2$ be a separable Hilbert space with
projectors $\Pi_i\colon \CH\to\CH_i$.  Let $(X_1,X_2)$ be an $\CH$-valued
Gaussian random variable with mean $m=(m_1,m_2)$ and positive
definite covariance operator~$\CC$ and define $\CC_{ij} = \Pi_i \CC \Pi_j^*$.
Then the conditional distribution of $X_1$ given~$X_2$ is Gaussian with mean
\begin{equ}[e:condmean]
  m_{1|2} = m_1 + \CC_{12}\CC_{22}^{-1} \bigl(X_2 - m_2\bigr)
\end{equ}
and covariance operator
\begin{equ}[e:condcov]
  \CC_{1|2} = \CC_{11} - \CC_{12}\CC_{22}^{-1}\CC_{21}.
\end{equ}
\end{lemma}

\begin{proof}
Note that by Lemma~\ref{L:covariance} the operator $\CC$ is trace
class.  Thus $\CC_{11}$ and $\CC_{22}$ are also trace class.
Let $\mu$ be the law of $X_2$ and let $\CH_0$ be the range of
$\CC_{22}^{1/2}$ equipped with the inner product
\begin{equ}
  \scal{x,y}_0 = \scal{C_{22}^{-1/2}x,C_{22}^{-1/2}y}.
\end{equ}
If we embed $\CH_0\hookrightarrow \CH_2$ via the trivial injection
$i(f)=f$, then we find $i^*(f) = C_{22}f$.  Since $i\circ i^* =
C_{22}$ is the covariance operator of~$\mu$, the space $\CH_0$ is its
reproducing kernel Hilbert space.  From \lem{hilbert-schmidt} we know
that $\CC_{12}\CC_{22}^{-1/2}$ is Hilbert-Schmidt from $\CH_2$ to
$\CH_1$ and hence bounded.  Thus we can define
\begin{equ}
  A = \CC_{12}\CC_{22}^{-1/2}\,\CC_{22}^{-1/2} = \CC_{12}\CC_{22}^{-1}
\end{equ}
as a bounded operator from~$\CH_0$ to~$\CH_1$.

Let $(\phi_n)_n$ be an orthonormal basis of~$\CH_2$.  Then $\psi_n =
C_{22}^{1/2}\phi_n$ defines an orthonormal basis on~$\CH_0$ and we get
\begin{equ}
  \sum_{n\in\N} \| A\psi_n \|^2_{\CH_1}
    = \sum_{n\in\N} \| \CC_{12}\CC_{22}^{-1}\,\CC_{22}^{1/2}\phi_n \|^2_{\CH_1}
    = \sum_{n\in\N} \| \CC_{12}\CC_{22}^{-1/2}\phi_n \|^2_{\CH_1}
    < \infty,
\end{equ}
where the last inequality comes from \lem{hilbert-schmidt}.  This
shows that the operator~$A$ is Hilbert-Schmidt on the reproducing
kernel Hilbert space~$\CH_0$.  Theorem~II.3.3
of~\cite{Dalecky-Fomin91} shows that $A$ can be extended in a
measurable way to a subset of~$\CH_2$ which has full measure, so that
\eref{e:condmean} is well-defined.

Now consider the process $Y$ defined by
\begin{equ}
\begin{pmatrix} Y_1 \\ Y_2 \end{pmatrix}
  = \begin{pmatrix}
      I_{\CH_1} & - A \\
      0_{\CH_2} & I_{\CH_2}
    \end{pmatrix}
    \begin{pmatrix} X_1 \\ X_2 \end{pmatrix}.
\end{equ}
This process is also Gaussian, but with mean
\begin{equ}
m^Y
  = \begin{pmatrix}
      I_{\CH_1} & - A \\
      0_{\CH_2} & I_{\CH_2}
    \end{pmatrix}
    \begin{pmatrix} m_1 \\ m_2 \end{pmatrix}
  = \begin{pmatrix} m_1 - A m_2 \\ m_2 \end{pmatrix}
\end{equ}
and covariance operator
\begin{equ}
\CC^Y
  = \begin{pmatrix}
      I_{\CH_1} & - A \\
      0_{\CH_2} & I_{\CH_2}
    \end{pmatrix}
    \begin{pmatrix}
      \CC_{11} & \CC_{12} \\
      \CC_{21} & \CC_{22}
    \end{pmatrix}
    \begin{pmatrix}
      I_{\CH_1} & 0_{\CH_2} \\
      - A^* & I_{\CH_2}
    \end{pmatrix}
  = \begin{pmatrix}
      \CC_{11}- \CC_{12} \CC_{22}^{-1} \CC_{21} & 0 \\
      0 & \CC_{22}
    \end{pmatrix}.
\end{equ}
This shows that $Y_1 = X_1 - \CC_{12}\CC_{22}^{-1}X_2$ and $Y_2=X_2$ are
uncorrelated and thus independent.  So we get
\begin{equs}
\E\bigl(X_1\bigm|X_2\bigr)
 &= \E\bigl(X_1 - \CC_{12}\CC_{22}^{-1}X_2\bigm|X_2\bigr)
    + \E\bigl(\CC_{12}\CC_{22}^{-1}X_2\bigm|X_2\bigr) \\
 &= \E\bigl(X_1 - \CC_{12}\CC_{22}^{-1}X_2\bigr)
    + \CC_{12}\CC_{22}^{-1}X_2 \\
 &= m_1 - \CC_{12}\CC_{22}^{-1}m_2 + \CC_{12}\CC_{22}^{-1}X_2.
\end{equs}
This proves~\eref{e:condmean} and a similar calculation gives
equality~\eref{e:condcov}.
\end{proof}

\begin{remark}\label{rem:just}
If we define as above $\CL = (-\CC)^{-1}$ and formally define
$\CL_{ij} = \Pi_i \CL \Pi_j^*$ (note that without additional
information on the domain of $\CL$ these operators may not be densely
defined), then a simple formal calculation shows that $m_{1|2}$ and
$\CC_{1|2}$ are expected to be given by
\begin{equ}[e:relcond]
m_{1|2} = m_1 - \CL_{11}^{-1} \CL_{12} \bigl(X_2 - m_2\bigr),\qquad
\CC_{1|2} = -\CL_{11}^{-1}.
\end{equ}
We now justify these relations in a particular situation which is
adapted to the case that will be considered in the remaining part of
this section.
\end{remark}

\begin{lemma}\label{lem:abstr}
  Consider the setup of \lem{gausscond} and \rem{rem:just} and
  assume furthermore that the following properties are satisfied:
\begin{claim}
\item[a.] The operator $\CL$ can be extended to a closed operator
  $\tilde \CL$ on $\Pi_1 \CD(\CL) \oplus \Pi_2 \CD(\CL)$.
\item[b.] Define the operators $\CL_{ij} = \Pi_i \tilde \CL \Pi_j^*$.
  Then, the operator $\CL_{11}$ is self-adjoint and one has $\ker
  \CL_{11} = \{0\}$.
\item[c.] The operator $-\CL_{11}^{-1} \CL_{12}$ can be extended to a
  bounded operator from $\CH_2$ into $\CH_1$.
\end{claim}
Then $\CC_{12} \CC_{22}^{-1}$ can be extended to a bounded operator
from $\CH_2$ into $\CH_1$ and one has $\CC_{12} \CC_{22}^{-1} =
-\CL_{11}^{-1} \CL_{12}$.  Furthermore, $\CC_{21}$ maps $\CH_1$ into
the range of $\CC_{22}$ and one has
\begin{equ}
\CL_{11}^{-1} x =  \bigl(\CC_{11} - \CC_{12} \CC_{22}^{-1} \CC_{21}\bigr) x,
\end{equ}
for every $x \in \CH_1$.
\end{lemma}

\begin{proof}
We first show that $\CC_{12} \CC_{22}^{-1} = -\CL_{11}^{-1} \CL_{12}$.
By property \textit{a.}\ and the definition of $\CL$, we have the
equality
\begin{equ}[e:prop1]
\tilde \CL \Pi_1^*\Pi_1 \CC x + \tilde \CL \Pi_2^*\Pi_2 \CC x = -x
\end{equ}
for every $x \in \CH$, and thus $\CL_{11} \CC_{12} x = - \CL_{12}
\CC_{22} x$ for every $x \in \CH_2$.  It follows immediately that
$\CL_{11} \CC_{12} \CC_{22}^{-1} x = -\CL_{12} x$ for every $x \in
\CR(\CC_{22})$. Since $\CR(\CC_{22})$ is dense in $\CH_2$, the
statement follows from assumptions \textit{b.} and~\textit{c.}

Let us now turn to the second equality.  By property~\textit{a.}\ the
operator $\CC_{21}$ maps $\CH_1$ into the domain of $\CL_{12}$ so that
\begin{equ}[e:prop2]
x = x - \CL_{12} \CC_{21} x + \CL_{12} \CC_{21} x
  = \CL_{11} \CC_{11} x + \CL_{12} \CC_{21} x,
\end{equ}
for every $x \in \CH_1$ (the second equality follows from an argument
similar to the one that yields \eref{e:prop1}).  Since the operator
$\CC_{22}^{-1}$ is self-adjoint, we know from \cite[p.~195]{Yos95FA}
that $(\CC_{12} \CC_{22}^{-1})^* = \CC_{22}^{-1} \CC_{21}$.  Since the
left hand side operator is densely defined and bounded, its adjoint is
defined on all of $\CH_1$, so that $\CC_{21}$ maps $\CH_1$ into the
range of $\CC_{22}$. It follows from \eref{e:prop2} that
\begin{equ}
x = \CL_{11} \CC_{11} x + \CL_{12} \CC_{22} \CC_{22}^{-1} \CC_{21} x,
\end{equ}
for every $x \in \CH_1$.  Using \eref{e:prop1}, this yields $x =
\CL_{11} \CC_{11} x - \CL_{11} \CC_{12} \CC_{22}^{-1} \CC_{21} x$, so
that $\CL_{11}^{-1}$ is an extension of $\CC_{11} - \CC_{12}
\CC_{22}^{-1} \CC_{21}$.  Since both of these operators are
self-adjoint, they must agree.
\end{proof}

\begin{corollary}
  Let $(X,Y)$ be Gaussian with covariance~$\CC$ and mean~$m$ on a
  separable Hilbert space $\CH = \CH_1 \oplus \CH_2$.  Assume
  furthermore that $\CC$ satisfies the assumptions of Lemmas
  \ref{gausscond} and~\ref{lem:abstr}.  Then, the conditional law of
  $X$ given $Y$ is given by the invariant measure of the ergodic SPDE
\begin{equ}[e:equcond]
\frac{dx}{dt}
  = \CL_{11} x - \CL_{11} \Pi_1 m
                + \CL_{12} \bigl(Y - \Pi_2 m\bigr) + \sqrt 2 \, \frac{dw}{dt},
\end{equ}
where $w$ is a cylindrical Wiener process on $\CH_1$ and the operators
$\CL_{ij}$ are defined as in~\lem{lem:abstr}.  SPDE~\eref{e:equcond}
is again interpreted in the mild sense~\eref{e:defmild}.
\end{corollary}

\begin{proof}
Note that $\CL_{11}^{-1} \CL_{12}$ can be extended to a bounded
operator by assumption and the mild interpretation
of~\eref{e:equcond} is
\begin{equ}[e:equmild]
x_t = M + e^{\CL_{11} t} (x_0 - M) + \sqrt 2 \int_0^t e^{\CL_{11}(t-s)}\,dw(s),
\end{equ}
with $M = \Pi_1 m - \CL_{11}^{-1} \CL_{12} \bigl(Y - \Pi_2 m\bigr)$.
The result follows by combining \lem{gausscond} and \lem{lem:abstr}
with~\lem{convergence}.
\end{proof}

These abstract results enable us to prove the main result of this
section.

\begin{proof}[of Theorem~\ref{kalman-SPDE}]
Consider a solution $(X,Y)$ to the SDE~\eref{e:kSDE}.
Introducing the shorthand notations
\begin{equ}
\Sigma_1 = (B_{11}B_{11}^*)^{-1},\qquad \Sigma_2 = (B_{22}B_{22}^*)^{-1},
\end{equ}
it follows by the techniques used in the proof of
Theorem~\ref{samplingLn} that the operator~$\CL$ corresponding to its
covariance is formally given by
\begin{equs}
\lhs
\begin{pmatrix}
  L_{11} & L_{12} \\
  L_{21} & L_{22}
\end{pmatrix}
 : = \begin{pmatrix}
      \d_u + A_{11}^* & A_{21}^* \\
      0 & \d_u
    \end{pmatrix}
    \begin{pmatrix}
      B_{11}B_{11}^*  & 0 \\
      0 & B_{22}B_{22}^*
    \end{pmatrix}^{-1}
    \begin{pmatrix}
      \d_u - A_{11} & 0 \\
      -A_{21} & \d_u
    \end{pmatrix} \\
  &= \begin{pmatrix}
      (\d_u + A_{11}^*) \Sigma_1 (\d_u - A_{11})
      - A_{21}^*\Sigma_2 A_{21}
    & A_{21}^*\Sigma_2 \d_u \\
      - \d_u \Sigma_2 A_{21}
    & \d_u \Sigma_2 \d_u
     \end{pmatrix}.
\end{equs}
In order to identify its domain, we consider~\eref{e:c3SPDEb} with
\begin{equ}
\Sigma =
  \begin{pmatrix}
    \Lambda & 0 \\
    0 & \Gamma
  \end{pmatrix}
\end{equ}
and we take the limit $\Gamma\to 0$.  This leads to the boundary
conditions
\minilab{e:k2SPDE}
\begin{equs}[e:k2SPDEb1][2]
\d_u x (0) &= A_{11}x(0) + (\Lambda\Sigma_1)^{-1} (x(0)-x^-),
  &\quad \d_u x (1) &= A_{11} x(1), \\
y(0) &= 0,
  &\quad \d_u y (1) &= A_{21} x(1).
\end{equs}
The domain of~$\CL$ is thus $H^2([0,1],\R^m \times \R^n)$, equipped
with the the homogeneous version of these boundary
conditions.

We now check that the conditions of \lem{lem:abstr} hold.  Condition
\textit{a.}\ is readily verified, the operator~$\tilde\CL$ being
equipped with the boundary conditions
\minilab{e:k2SPDE}
\begin{equs}[2][e:k2SPDEb2]
\d_u x (0) &= A_{11}x(0) + (\Lambda\Sigma_1)^{-1} x(0),
  &\quad \d_u x (1) &= A_{11} x(1), \\
y(0) &= 0,&\quad \Pi \d_u y(1) &= 0,
\end{equs}
where $\Pi$ is the projection on the orthogonal complement of
the range of $A_{21}$.
Note that the operator $\tilde \CL$ is closed, but no longer
self-adjoint (unless $A_{21} = 0$).  The operator $\CL_{11}$ is therefore
given by
\begin{equ}
\CL_{11} = \bigl(\d_u  + A_{11}^*\bigr) \Sigma_1 \bigl(\d_u - A_{11}\bigr)
      - A_{21}^*\Sigma_2 A_{21},
\end{equ}
equipped with the boundary condition
\begin{equ}
\d_u x (0) = A_{11}x(0) + (\Lambda\Sigma_1)^{-1} x(0),
  \quad \d_u x (1) = A_{11} x(1).
\end{equ}
It is clear that this operator is self-adjoint.  The fact that its
spectrum is bounded away from $0$ follows from the fact that the form
domain of~$\CL$ contains $\Pi_1^* \Pi_1 \CD(\CL)$ and that there is a
$c>0$ with $\scal{a, \CL a} \leq - c \|a\|^2$ for all $a \in
\CD(\CL)$.  Thus condition~\textit{b.}\ holds.

The operator $\CL_{12}$ is given by the first-order
differential operator $A_{21}^*\Sigma_2 \d_u$ whose domain is given by
functions with square-integrable second derivative that vanish at $0$.
Since the kernel of $\CL_{11}^{-1}$ has a square-integrable
derivative, it is easy to check that $\CL_{11}^{-1}\CL_{12}$ extends
to a bounded operator on $\CH$, so
that condition \textit{c.}\ is also verified.

We can therefore apply \lem{lem:abstr} and \lem{convergence}. The
formulation of the equation with inhomogeneous boundary conditions is
an immediate consequence of \rem{rem:BC}: a short calculation to
remove the inhomogeneity in the boundary conditions~\eref{e:kSPDEb} and
change the inhomogeneity in the PDE~\eref{e:kSPDEa} shows
that~\eref{e:kSPDE} can be written in the form~\eref{e:equcond}
or~\eref{e:equmild} with the desired value for~$M$, the conditional
mean.  Since~$\CL_{11}$ is indeed the conditional covariance operator,
the proof is complete.
\end{proof}

\begin{remark}
  For $Y$ solving~\eref{e:kSDE} the derivative $\frac{dY}{du}$ only
  exists in a distributional sense (it is in the Sobolev space
  $H^{-1/2-\eps}$ for every $\eps>0$).  But the
  definition~\eref{e:defmild} of a mild solution which we use here
  applies the inverse of the second order differential
  operator~$\CL_{11}$ to~$\frac{dY}{du}$, resulting in an element of
  $H^{3/2-\eps}$ in the solution.
\end{remark}

\begin{remark}
  Denote by $x(t,u)$ a solution of the SPDE~\eref{e:kSPDE} and write
  the mean as $\bar x(t,u) = \E x(t,u)$.  Then, as $t \to \infty$,
  $\bar x(t,u)$ converges to its limit $\tilde x(u)$ strongly in
  $L^2([0,1],\R^m)$ and $\tilde x(u)$ must coincide with the
  Kalman-Bucy filter/smoother.  This follows from the fact that
  $\tilde x$ equals $\E (X\,|\,Y)$.  It is instructive to demonstrate
  this result directly and so we do so.

The mean~$\tilde x(u)$ of the invariant measure of~\eref{e:kSPDE}
satisfies the linear two point boundary value problem
\begin{subequations} \label{e:kBVP}
\begin{gather}
\begin{split}
&(\frac{d}{du} + A_{11}^*)
      (B_{11}B_{11}^*)^{-1}(\frac{d}{du} - A_{11}) \tilde x(u) \\
&\hskip5mm
 + A_{21}^*(B_{22}B_{22}^*)^{-1}
 \Bigl( \frac{dY}{du} - A_{21} \tilde x(u)\Bigr) = 0 \qquad \forall u \in (0,1),
\end{split} \label{e:kBVPa} \\
\frac{d}{du}\tilde x(0)
  = A_{11}\tilde x(0) + B_{11}B_{11}^* \Lambda^{-1} \bigl(\tilde x(0)-x^-\bigr),
                                                    \label{e:kBVPb} \\
\frac{d}{du}\tilde x(1) = A_{11}\tilde x(1). \label{e:kBVPc}
\end{gather}
\end{subequations}

The standard implementation of the Kalman filter is to calculate the
conditional expectation $\hat X(u) = \E\bigl(X(u) \bigm| Y(v), 0\leq
v\leq u\bigr)$ by solving the initial value problem
\begin{equs}
\frac{d}{du}S(u)
  &= A_{11} S(u)+S(u)A_{11}^* - S(u) A_{21}^*(B_{22}B_{22}^*)^{-1}A_{21} S(u)
       + B_{11} B_{11}^* \\
S(0)
  &= \Lambda \label{e:volat}
\end{equs}
and
\begin{equs}
\frac{d}{du} \hat X(u)
  &= \bigl( A_{11}-S(u) A_{21}^*(B_{22}B_{22}^*)^{-1}A_{21} \bigr) \hat X
       + S(u) A_{21}^*(B_{22}B_{22}^*)^{-1} \frac{dY}{du} \\
\hat X(0)
  &= x^-.\label{e:filter}
\end{equs}
The Kalman smoother~$\tilde X$, designed to find
$\tilde X(u) = \E\bigl(X(u) \bigm| Y(v),  0\leq v\leq 1\bigr)$,
is then given by the backward sweep
\begin{equs}
\frac{d}{du}\tilde X(u)
  &= A_{11} \tilde X(u) + B_{11} B_{11}^* S(u)^{-1}
                               \bigl(\tilde X(u)-\hat X(u)\bigr)
        \qquad\forall u\in(0,1) \\
\tilde X(1) &= \hat X(1). \label{e:smoother}
\end{equs}
See~\cite[section~6.3 and exercise~6.6]{Oksendal98} for a reference.
We wish to demonstrate that $\tilde x(u)=\tilde X(u)$.

Equation~\eref{e:smoother} evaluated for $u=1$ gives
equation~\eref{e:kBVPc}.  When evaluating~\eref{e:smoother} at~$u=0$
we can use the boundary conditions from \eref{e:volat}
and~\eref{e:filter} to get equation~\eref{e:kBVPb}.  Thus it remains
to show that $\tilde X(u)$ satisfies equation~\eref{e:kBVPa}.  We
proceed as follows: equation~\eref{e:smoother} gives
\begin{equs}
\lhs (\frac{d}{du} + A_{11}^*) (B_{11}B_{11}^*)^{-1}
              (\frac{d}{du} - A_{11}) \tilde X \\
  &= (\frac{d}{du} + A_{11}^*) (B_{11}B_{11}^*)^{-1}
              B_{11} B_{11}^* S^{-1} \bigl(\tilde X-\hat X\bigr) \\
  &= (\frac{d}{du} + A_{11}^*)
              S^{-1} \bigl(\tilde X-\hat X\bigr)
\end{equs}
and so
\begin{equs}[eq:n1][1]
\lhs \bigl(\frac{d}{du} + A_{11}^*\bigr) (B_{11}B_{11}^*)^{-1}
              (\frac{d}{du} - A_{11}) \tilde X \\
 &= \Bigl(A_{11}^*S^{-1} + \frac{d}{du} S^{-1}\Bigr) \bigl(\tilde X-\hat X\bigr)
     + S^{-1} \frac{d}{du} \bigl(\tilde X-\hat X\bigr).
\end{equs}
We have
\begin{equation*}
\frac{d}{du} S^{-1}
  = - S^{-1} \frac{dS}{du} S^{-1}
\end{equation*}
and hence, using equation~\eref{e:volat}, we get
\begin{equation}
\label{eq:n2}
\frac{d}{du} S^{-1}
  = - S^{-1} A_{11} - A_{11}^* S^{-1} + A_{21}^*(B_{22}B_{22}^*)^{-1}A_{21}
     - S^{-1} B_{11} B_{11}^* S^{-1}.
\end{equation}
Subtracting \eref{e:filter} from~\eref{e:smoother} leads to
\begin{align}
\label{eq:n3}
S^{-1} \frac{d}{du} \bigl(\tilde X-\hat X\bigr)
  = S^{-1} A_{11} \bigl(\tilde X-\hat X\bigr)
    &+ S^{-1} B_{11} B_{11}^* S^{-1} \bigl(\tilde X-\hat X\bigr)
\notag\\
  &- A_{21}^*(B_{22}B_{22}^*)^{-1}
                \Bigl( \frac{dY}{du} - A_{21} \hat X \Bigr).
\end{align}
By substituting \eqref{eq:n2}, \eqref{eq:n3} into
\eqref{eq:n1} and collecting all the terms we find
\begin{equ}
(\frac{d}{du} + A_{11}^*) (B_{11}B_{11}^*)^{-1}
              (\frac{d}{du} - A_{11}) \tilde X
  = - A_{21}^*(B_{22}B_{22}^*)^{-1}
                \Bigl( \frac{dY}{du} - A_{21} \tilde X \Bigr)
\end{equ}
which is equation~\eref{e:kBVPa}.

We note in passing that equations \eref{e:volat} to~\eref{e:smoother}
constitute a factorisation of the two-point boundary value
problem~\eref{e:kBVP} reminiscent of a continuous LU-factorisation
of~$\CL_{11}$.
\end{remark}


\section{Numerical Approximation of the SPDEs and Sampling}
\label{S:disc}

A primary objective when introducing SPDEs in this paper, and in the
nonlinear companion~\cite{Hairer-Stuart-VossII}, is to construct MCMC
methods to sample conditioned diffusions.  In this section we
illustrate briefly how this can be implemented.

\medskip

If we discretise the SDE~\eqref{e:L2SDE} in time by the
$\theta$-method, we obtain the following implicitly defined mapping
from $(x^k,\xi^k)$ to $x^{\star}$:
\begin{equation*}
\frac{x^{\star}-x^k}{\Delta t}=
\Bigl(\theta \CL x^{*}+(1-\theta)\CL x^k \Bigr)-\CL m
+\sqrt\frac{2}{\Delta t}\xi^k,
\end{equation*}
where $\xi^k$ is a sequence of i.i.d Gaussian random variables  in 
$\CH$ with covariance operator $I$ (i.e. white noise in in $\CH$).
The Markov chain implied by the map is well-defined on $\CH$ 
for every $\theta \in [\frac12,1]$.
This Markov chain can be used as a proposal distribution for
an MCMC method, using the Metropolis-Hastings criterion to
accept or reject steps. To make a practical algorithm it is
necessary to discretise in the Hilbert space $\CH$, as well as
in time $t$. This idea extends to nonlinear problems.

Straightforward calculation using the Karhunen-Lo\`{e}ve expansion, similar
to the calculations following Lemma \ref{convergence},  shows that the 
invariant measure of the SPDE \eqref{e:L2SDE} is preserved if the SPDE 
is replaced by
\begin{equ}[e:L3SDE]
  \frac{dx}{dt} = -x + m + \sqrt{2\CC} \,\frac{dw}{dt}.
\end{equ}
Such pre-conditioning of Langevin equations can be beneficial
algorithmically because it equalises convergence rates
in different modes. This in turn allows for optimisation
of the time-step choice for a Metropolis-Hastings algorithm
across all modes simultaneously. We illustrate this issue
for the linear Gaussian processes of interest here.

Equation \eqref{e:L3SDE} can be discretised in time by the 
$\theta$-method to obtain the following implicitly defined mapping
from $(x^k,\xi^k)$ into $x^{\star}$: 
\begin{equation*}
\frac{x^{\star}-x^k}{\Delta t}=
-\Bigl(\theta x^{\star}+(1-\theta)x^k \Bigr)+m
+\sqrt\frac{2}{\Delta t}\xi^k.
\end{equation*}
Now $\xi^k$ is a sequence of i.i.d.\ Gaussian random variables in
$\CH$ with covariance operator~$\CC$.  Again, this leads to a
well-defined Markov chain on~$\CH$ for every $\theta \in [\frac12,1]$.
Furthermore the invariant measure is $\CC/(1+(\theta-\frac12)\Delta
t)$.  Thus the choice $\theta=\frac12$ has a particular advantage: it
preserves the exact invariant measure, for all $\Delta t >0$.  (These
observations can be justified by using the Karhunen-Lo\`eve
expansion).  Note that
$$\Bigl(1+\theta \Delta t\Bigr)x^{\star}=
\Bigl(1-(1-\theta)\Delta t\Bigr)x^k+\sqrt{2\Delta t}\xi_k.$$
When $\theta=\frac12$, choosing $\Delta t=2$ generates
independent random variables which therefore sample the invariant
measure independently.  This illustrates in a simple Gaussian setting
the fact that it is possible to choose a globally optimal time-step
for the MCMC method.  To make a practical algorithm it is necessary to
discretise in the Hilbert space $\CH$, as well as in time $t$. The
ideas provide useful insight into nonlinear problems.


\section{Conclusions}
\label{S:conc}

In this text we derived and exploited a method to construct linear
SPDEs which have a prescribed Gaussian measure as their stationary
distribution.  The fundamental relation between the diffusion
operator~$\CL$ in the SPDE and the covariance operator~$\CC$ of the
Gaussian measure is $\CL = (-\CC)^{-1}$ and, using this, we showed
that the kernel of the covariance operator (the covariance function)
is the Green's functions for~$\CL$.  We illustrated this technique by
constructing SPDEs which sample from the distributions of linear SDEs
conditioned on several different types of observations.

These abstract Gaussian results were used to produce some interesting
results about the structure of the Kalman-Bucy filter/smoother.
Connections were also made between discretisations of the resulting
SPDEs and MCMC methods for the Gaussian processes of interest.

In the companion article~\cite{Hairer-Stuart-VossII} we build on the
present analysis to extend this technique beyond the linear case.
There we consider conditioned SDEs where the drift is a gradient (or
more generally a linear function plus a gradient).  The resulting
SPDEs can be derived from the SPDEs in the present text by the
addition of an extra drift term to account for the additional
gradient. The stationary distributions of the new nonlinear SPDEs are
identified by calculating their Radon-Nikodym derivative with respect to the
corresponding stationary distributions of the linear equations as
identified in the present article; this is achieved via the
Girsanov transformation.

The Girsanov transformation is used to study the connection between
SPDEs and bridge processes in \cite{Reznikoff-VandenEijnden05}; it is
also used to study Gibbs measures on $\R$ in~\cite{BL00}.  However the
results concerning bridges in this paper are not a linear subcase of
those papers because we consider non-symmetric drifts (which are hence
not gradient) and covariance of the noise which is not proportional to
the identity.  Furthermore the nonlinear results
in~\cite{Hairer-Stuart-VossII} include the results stated in
\cite{Reznikoff-VandenEijnden05} as a subset, both because of the form
of the nonlinearity and noise, and because of the wide-ranging forms
of conditioning that we consider.


\bibliographystyle{alpha}
\bibliography{gaussian}

\end{document}